\RequirePackage{fix-cm}
\documentclass[smallextended,natbib,
a4paper]{svjour3}       
\smartqed  
\usepackage{graphicx}
\usepackage{amssymb}
\usepackage{mathptmx}      
\usepackage{tabularx}
\usepackage{booktabs}
\usepackage[inline]{enumitem}
\usepackage[draft]{hyperref}
\usepackage{breakcites}

\renewcommand{\cite}{\citep}

\newcolumntype{Y}{>{\centering\arraybackslash}X}

\journalname{Axiomathes}

\title{Deep Disagreement in Mathematics}
\author{Andrew Aberdein}
\institute{A. Aberdein \at
School of Arts and Communication, Florida Institute of Technology,
150 West University Blvd,
Melbourne, FL 32901, USA \\
\email{aberdein@fit.edu}}

\begin{document}
\maketitle

\abstract{
Disagreements that resist rational resolution, often termed “deep disagreements”, have been the focus of much work in epistemology and informal logic. 
In this paper, I argue that they also deserve the attention of philosophers of mathematics. 
I link the question of whether there can be deep disagreements in mathematics to a more familiar debate over whether there can be revolutions in mathematics. I propose an affirmative answer to both questions, using the controversy over Shinichi Mochizuki’s work on the $abc$ conjecture as a potential example of both phenomena. I conclude by investigating the prospects for the resolution of mathematical deep disagreements in virtue-theoretic approaches to informal logic and mathematical practice.
}

\section{Mathematical Revolutions}\label{sec:Revolutions}

Perhaps the earliest 
``fully modern and explicit statement concerning a revolution in mathematics'' 
was made in 1720 by Bernard de Fontenelle, permanent secretary of the Acad\'emie Royale des Sciences \cite[90]{Cohen85}.
In his eulogy for Michel Rolle,
Fontenelle characterized the innovations of the infinitesimal calculus as ``une r\'evolution bien marqu\'ee'' \cite[quoted in][214]{Cohen85}.
Bernard Cohen identifies many subsequent mathematicians and philosophers who have spoken in terms that imply the occurrence of mathematical revolutions, including Immanuel Kant, 
William Kingdon Clifford, 
Georg Cantor,
G\"osta Mittag-Leffler,
Morris Kline, and 
Beno\^it Mandelbrot \citetext{\citealp[490 f.]{Cohen85}; see also \citealp[110]{Francois10}}.
However, as Cohen also notes, other distinguished authorities have come to a contrary conclusion. 
Indeed, Michael Crowe declared as a ``law'' that ``Revolutions never occur in mathematics'', 
citing Joseph Fourier, Hermann Hankel, and Clifford Truesdell in support of a perspective on mathematical knowledge as inherently cumulative \cite[165]{Crowe75}.
Raymond Wilder offered a slightly more nuanced version of Crowe's ``law'', that ``Revolutions may occur in the metaphysics, symbolism and methodology of mathematics, but not in the core of mathematics'' \cite[142]{Wilder81}.


Despite Crowe's initial disavowal of mathematical revolutions, 
he subsequently characterized the proposition that ``mathematics is cumulative'' as a ``misconception'' and 
eventually concluded that the ``question of whether revolutions occur in mathematics is in substantial measure definitional'' \citetext{\citealp[263]{Crowe88}; \citealp[316]{Crowe92}}. 
Crowe could 
thereby continue to deny the existence of mathematical revolutions, but only if
``a necessary characteristic of a revolution is that some previously existing entity (be it king, constitution, or theory) must be overthrown and irrevocably discarded'' \cite[165]{Crowe75}.
Such support as he later offered mathematical revolutions would remain subject to this qualification:
\begin{quote}
{Massive areas of mathematics have, for all practical purposes, been abandoned}. The nineteenth-century mathematicians who extended two millennia of research on conic section theory have now been forgotten; invariant theory, so popular in the nineteenth century, fell from favor. Of the hundreds of proofs of the Pythagorean theorem, nearly all are now nothing more than curiosities. In short, although many previous areas, proofs, and concepts in mathematics have persisted, others are now abandoned. Scattered over the landscape of the past of mathematics are {numerous citadels}, once proudly erected, but which, {although never attacked, are now left unoccupied} by active mathematicians \cite[263]{Crowe88}.
\end{quote}
In other words, there are 
mathematical revolutions, but they are revolutions of a 
special sort. 
In the 
natural sciences, 
revolution sweeps away the old concepts---phlogiston is gone 
and so much the worse for Priestley; 
Lavoisier 
invoked fundamentally different theoretical posits 
to make sense of 
his experiments. Mathematicians are seldom in that 
situation, or so it would seem. 
Instead, they move on. What had been an area of widespread fascination becomes 
deserted. Crowe's abandoned citadels 
were never attacked, 
whereas phlogiston was roundly attacked and irrevocably discarded. 
Nothing like that happened to conic sections, 
yet nobody has regarded them as an active area of research in 
150 years. 
Hence mathematical revolutions are unlike other revolutions---or so it is argued.

In my earlier treatment of revolutions in logic I framed this distinction 
as between glorious and inglorious revolutions, by allusion to the Glorious Revolution of 1688 \cite[618 f.]{Aberdein09a}. 
The defeat of the Jacobites and the ascent of William and Mary was 
consequential with respect to the balance of power, but Britain remained a monarchy with a parliament throughout. 
The crucial shift was in the relative significance of these different components, not the loss of 
a crucial component or the introduction of 
an unprecedented new component.%
\footnote{Donald Gillies draws a similar distinction using different historical examples: Franco-British revolutions in Gillies' terms are my glorious revolutions, whereas his Russian is my inglorious revolution \cite[5]{Gillies92}.}
I also proposed the concept of a paraglorious revolution: this adds something which is not commensurable with what went before. 
The paraglorious revolution 
thereby resembles an inglorious revolution 
with the arrow of time reversed. 

Another of the misconceptions that Crowe diagnoses in the historiography of mathematics is that ``the methodology of mathematics is radically different from the methodology of science'' \cite[271]{Crowe88}. If we accept a methodological continuity between mathematics and science we should also expect not only that both fields should exhibit revolutions but also that both sorts of revolution should be amenable to similar analysis.
The study of scientific revolutions has become irresistibly linked to the work of Thomas Kuhn, for whom they are paradigm shifts or transitions between incommensurable scientific theories. 
Kuhn says many 
things 
about incommensurability 
and fitting them together is a complicated task. 
Here is one 
significant distinction which Kuhn does not make entirely explicit: 
\begin{quote}
\begin{description}
\item[Methodological Incommensurability (MI)] There are no objective criteria of theory evaluation. The familiar criteria of evaluation, such as simplicity and fruitfulness, are not a fixed set of rules. Rather, they vary with the currently dominant paradigm. 
\item[Taxonomic Incommensurability (TI)] Periods of scientific change ({in particular, revolutionary change}) that exhibit TI are scientific developments in which existing concepts are {replaced} with new concepts that are {incompatible} with the older concepts. The new concepts are incompatible with the old concepts in the following sense: two competing scientific theories are conceptually incompatible (or incommensurable) {just in case they do not share the same ``lexical taxonomy''}. A lexical taxonomy contains the structures and vocabulary that are used to state a theory 
\cite[362, internal citations omitted]{Mizrahi14b}.
\end{description}
\end{quote}
How do these two understandings of incommensurability fit the mathematical case? 
Methodological incommensurability is perhaps the more familiar. 
Here, for example, is Yehuda Rav discussing 
the diversity of mathematical practices: 
``because of the historical and \emph{methodological} wealth of mathematical proof practices (plural), any attempt to encapsulate such multifarious practices in a unique and uniform one-block perspective is bound to be defective'' \cite[299, emphasis in original]{Rav07}. 
As examples, he notes the differences 
between ``(a) historical epochs; (b) different (and at times opposing) acceptability criteria of logical reasoning in mathematics; and (c) distinct admissibility criteria of types of proof in different branches of mathematics'' \citetext{ibid.}. 
This looks 
like methodological incommensurability: no objective criteria common across every domain. 
However, Rav is not claiming that any of these differences suffice for a revolution (although some revolutions may exhibit such methodological changes).
As we shall see, taxonomic incommensurability is a closer match to the discontinuity characteristic of mathematical revolutions.
For example, Ken Manders describes what he calls 
domain extension by existential closure and model completion, a process whereby a theory is supplemented with new elements, strictly incompatible with the existing elements, 
such as adding a square root of $-1$ to the real numbers, and thereby ensuring that all algebraic equations have solutions. As Manders explains,
``Model completion relates an established theoretical context to a new one, intimately connected: we retain individuals and certain of their relationships and properties in the original domain. But some issues about the original context become more understandable in the model-completed one'' \cite[557]{Manders89}.
Innovations of this sort transformed mathematics in the nineteenth century 
(\citealp{Gray92}; see also \citealp{Cantu13,Bellomo21}).

To give these ideas a little more precision let us make some definitions. We shall assume without further analysis that we are comparing theories which are comprised of multiple components (concepts, ideas, techniques,\dots), some of which are of 
key significance to the practitioners of that theory; that is, their absence would alter the theory in some fundamental way, at least in the eyes of its practitioners. We may then define some relationships between sets of key components and thereby between the respective theories:
\begin{definition}[Isomorphism]
Two sets of key components $K_i$ and $K_j$ are \emph{isomorphic} 
($K_{i}\cong K_{j}$) when there exists a 
correspondence between their elements that preserves meaning for the practitioners of the theories of which $K_i$ and $K_j$ are a part.
\end{definition}
\begin{definition}[Isomorphic embedding]
$K_i$ 
\emph{embeds isomorphically} in $K_j$ ($K_{i}\hookrightarrow K_{j}$) if and only if there exists some proper subset of $K_j$, $K_{j}^{\prime}$, such that $K_{i}\cong K_{j}^\prime$. 
\end{definition}
With these definitions in place, we may now give a somewhat more precise account of the distinction drawn above between different types of revolution:
\begin{definition}[Revolutions]
Where $K_{i}$ are the key components of theory $T_{i}$, succession from theory $T_{n}$ to theory $T_{n+1}$ is 
a \emph{null revolution} when $K_{n}= K_{n+1}$;
a \emph{glorious revolution} when $K_{n}\neq K_{n+1}$ and $K_{n}\cong K_{n+1}$;
a \emph{paraglorious revolution} when $K_{n}\ncong K_{n+1}$ and $K_{n}\hookrightarrow K_{n+1}$; and
an \emph{inglorious revolution} when $K_{n}\ncong K_{n+1}$ and $K_{n}\not\hookrightarrow K_{n+1}$. 
\end{definition}
As a special case of inglorious revolution, we can consider revolutions which simply restrict the set of key terms, such that $K_{n}\hookleftarrow K_{n+1}$, and are thereby dual to paraglorious revolutions. More generally, inglorious revolutions both add and remove key terms, such that there is some common set $K_{n}^{\prime}$, not necessarily itself the subject of a seriously defended theory, such that $K_{n}\hookleftarrow K_{n}^{\prime}\hookrightarrow K_{n+1}$.

These definitions raise several additional questions, which I have pursued at greater length elsewhere
\citetext{\citealp[619 ff.]{Aberdein09a}; \citealp[136 f.]{Aberdein17}}.
Firstly, 
the definitions turn on 
what the practitioners of the theories at issue 
take to be key and what they consider to have been preserved in moving from one theory to another.
A full accounting of these matters for any such transition is potentially a delicate and 
substantial task. 
Fortunately, 
there is often room for consensus between adherents of pre- and post-revolutionary theories: 
Priestley and Lavoisier could presumably agree that phlogiston was key to Priestley's theory and not preserved in Lavoisier's.
Further complicating this question is the issue of lexical drift:
the meaning of persistent terminology can change beyond recognition and new terminology can be used to refer to old concepts. 
Thus 
lexical taxonomy cannot be merely lexical; it must accommodate such potentially distracting shifts in the lexicon.
Detailed historical studies must be sensitive to all these factors \cite[see, for example,][]{Barany18}.
A further question concerns the status of a sequence of consecutive revolutions: if they are all of the same type, can they be treated as a single revolution of that type? Clearly not in general, as a pair of inglorious revolutions might in theory cancel each other out. Moreover, some sequences of individually glorious (or paraglorious) revolutions might be plausibly treated as collectively inglorious, for example if a sorites-like series of small changes in the meaning of key components adds up to a bigger change.

As we have seen, many mathematicians, and philosophers or historians of mathematics, have denied that there can be revolutions in mathematics. Some, such as (the later) Crowe and Gillies, have admitted what I have defined as glorious revolutions. My contention is that inglorious and paraglorious revolutions are also to be found.
What examples of such mathematical revolutions are there?
One problem for the identification of revolutions in mathematics is that of scale. While revolutions of Copernican magnitude may be found in many natural sciences, mathematical revolutions with such sweeping effects are harder to identify. This has been taken---erroneously, I would argue---as evidence that mathematics lacks true revolutions.
In the first place, this criticism rests on a 
narrow reading of Kuhn's characterization of revolution, 
as implicitly excluding ``small-scale `micro-revolutions'\,'' \cite[47]{Toulmin70}.
Without engaging the vexed questions of Kuhnian exegesis, I wish to endorse his view that such shifts are indeed revolutions:
``a little studied type of conceptual change which occurs \emph{frequently} in science and is fundamental to its advance'' \cite[249 f., my emphasis]{Kuhn70}.
Kuhn later gave an indication of the scale of such conceptual changes in defining a paradigm as ``what the members of a scientific community and they alone share'', 
where such communities may comprise ``perhaps 100 members, sometimes significantly fewer'' \cite[460; 462]{Kuhn77}.
Mathematical research communities are typically at the low end of this range, but certainly within it:
``a few dozen (at most a few hundred)
'' \cite[35]{Davis80}.%
\footnote{The activity of such communities is, of course, the ``normal science'' which Kuhn contrasts with scientific revolutions. Understanding the normal science of mathematical research communities is a key task of the philosophy of mathematical practice, but complementary to the focus of the present paper.
For one specific such project, see \cite{Petersen14}.} 
Secondly, there have been large-scale mathematical revolutions.
Bruce Pourciau has argued that Brouwerian intuitionism, which, if adopted, would have required wholesale revision of results treated as certain by prior mathematicians, is a (failed) Kuhnian revolution \cite{Pourciau00}.
In my earlier treatment of mathematical revolutions I proposed several 
further putative cases of mathematical revolution \cite[140 ff.]{Aberdein17}.
The shift in antiquity from rational to real numbers---the original case of incommensurability---
is at least a paraglorious revolution and perhaps also inglorious. 
Other examples of substantial mathematical revolutions proposed by historians include
non-Euclidean geometry \cite{Ashkenazi14} and
the origins of algebraic geometry \cite{Oaks18}.

Examples of smaller scale revolutions arise from 
the wide-ranging, so-called architectural conjectures that structure much of contemporary mathematics \cite[see][]{Mazur97}. 
When a whole research programme is based on a conjecture that turns out to be false, 
then the eventual failure of that conjecture must be a strictly inglorious revolution.
Such collapsing conjectures are proofs by contradiction writ large; they describe ``things that seem to exist but that, when we study them deeply enough, turn out not to exist after all'' \cite[18]{Propp19}.
In such cases, the failure of the conjecture requires the exclusion of elements of the set of key components.%
\footnote{For example, consider the collapse of the ``world without end hypothesis'', which posited the existence of a certain type of map, termed $\theta_j$, between higher-dimensional spheres 
for all dimensions $2^{j+1}-2$. The proof that no such maps exist for $j\geq7$ demolished an elaborate structure of conjectures comprising a theory for which the $\theta_j$ for $j\geq7$ were key components. For further discussion, see \cite[143 ff.]{Aberdein17}.} 
Conversely, shifts resulting from a rapid advance in key conceptual content 
are 
paraglorious revolutions (if they succeed; if they do not they are also failed architectural conjectures). 
Lastly, the shift from modern to contemporary mathematics has involved numerous conceptual innovations, such as transference, reflection, and gluing \cite[27]{Zalamea12}. 
Each of these steps proceeds by strictly augmenting the store of key concepts and techniques, and may therefore be seen as paraglorious revolutions.
Moreover, as the older concepts and techniques fall into desuetude in the manner of Crowe’s abandoned citadels, they may be expected to steadily lose their key status. 
Thus, taken collectively, the transformation of mathematics over the last couple of centuries may be seen as 
a large-scale revolution that is locally paraglorious and globally inglorious.

\section{Deep Disagreement}\label{sec:DD}
What is deep disagreement?
As a term of art, it originates with Robert Fogelin \citeyearpar{Fogelin85}. 
He stresses 
that deep disagreement is 
an epistemic consideration, 
unrelated to strength of feeling. 
The resistance to resolution is not just a contingent feature of the individual circumstances of the dispute:
\begin{quote}
A disagreement can be {intense without being deep}. A disagreement can also be {unresolvable without being deep}. I can argue myself blue in the face trying to convince you of something without succeeding. The explanation might be that one of us is {dense or pig-headed}. And this is a matter that could be established beyond doubt to, say, an impartial spectator. But we get a very different sort of disagreement when it proceeds from {a clash in underlying principles}. Under these circumstances, the parties may be unbiased, free of prejudice, consistent, coherent, precise and rigorous, yet still disagree. And disagree profoundly, not just marginally 
\cite[5]{Fogelin85}.
\end{quote}
Fogelin invokes Putnam and Wittgenstein to link his ``underlying principles'' to 
their rules or framework propositions. 
Others have made ``hinge commitments'' the basis of contention \cite{Pritchard19}.
Michael Lynch has focussed instead on epistemic principles. He offers a fourfold set of criteria for 
when a disagreement is deep:
\begin{quote}
\begin{enumerate}
\item \emph{Commonality}: The parties to the disagreement {share} common epistemic {goal(s)}.
\item \emph{Competition}: If the parties affirm distinct principles with regard to a given domain, those principles 
\begin{enumerate}
\item pronounce {different methods} to be the most reliable in a given domain; and 
\item these methods are capable of producing {incompatible beliefs} about that domain.
\end{enumerate}
\item \emph{Non-arbitration}: There is {no further epistemic principle}, accepted by both parties, which would settle the disagreement.
\item \emph{Mutual Circularity}: The epistemic principle(s) in question can be justified only by means of an epistemically {circular argument} \cite[265]{Lynch10}.
\end{enumerate}
\end{quote}
These criteria presume that the disagreement is over an epistemic principle, in other words a \emph{direct} deep disagreement. But disagreements can also be \emph{indirectly} deep, when the parties ostensibly disagree over something 
other than an epistemic principle, but their disagreement commits them to also disagreeing directly (whether or not they appreciate this) \cite[988]{Ranalli19a}.
For a disagreement 
to be deep 
on Lynch's account, 
all four conditions must be met. 
If the parties are trying to do different things 
or are actually talking about different 
stuff, they are not in disagreement; 
and if there is a further epistemic principle they can appeal to to resolve the issue
or a non-question-begging 
justification for one or other epistemic principle, their disagreement is not deep.


In place of epistemic principles, framework propositions, or 
hinge commitments, it is 
also plausible to see deep disagreement in Kuhnian terms as a clash of incommensurable paradigms. 
Recall the distinction drawn in the previous section between methodological and taxonomic incommensurability.
On a methodological approach to incommensurability, any criteria you might use for theory appraisal will be relative to context, so there is no higher authority that can arbitrate disputes. 
On a taxonomic approach to incommensurability, 
paradigms are incompatible 
because they employ different vocabularies. 

So how do deep disagreements connect to the account of mathematical revolutions outlined in the previous section? 
My proposal is that 
agents that disagree with each other are in deep disagreement if they endorse theories with sets of key components that are not isomorphic. 
Maybe there are other ways of being in deep disagreement; I claim 
only that this is sufficient for a deep disagreement. 

\begin{definition}[
Deep Disagreement]\label{def:DD}
Agents $S_i$ and $S_j$ are in deep disagreement about some subject matter if they endorse theories $T_i$ and $T_j$ of that subject matter 
such that $K_{i}\ncong K_{j}$.
\end{definition}
Does this definition satisfy Lynch's criteria? Let us take them in turn.
\begin{enumerate*}[label=(\arabic*)]
\item By definition, $T_i$ and $T_j$ concern the same phenomena, satisfying Commonality.
\item  If $K_{i}\ncong K_{j}$, then $T_i$ and $T_j$ necessarily contain incompatible beliefs. Competition also requires that these incompatible beliefs be mediated by a difference of method, which could be taken to require that $T_i$ and $T_j$ are not just taxonomically incommensurable, but also methodologically incommensurable. These forms of incommensurability are conceptually independent (even though they mostly coincide). Thus Definition \ref{def:DD} could be strengthened to ensure that $T_i$ and $T_j$ are methodologically incommensurable too, thereby fully satisfying Competition. However, without seeking to settle this issue, I will employ the weaker form of Definition \ref{def:DD}.
\item Any further epistemic principle, common to both $K_i$ and $K_j$ which could arbitrate the dispute between them, would thereby demonstrate that $K_{i}\cong K_{j}$, hence Definition \ref{def:DD} implies Non-arbitration.
\item We would not expect every element of $K_i$ or $K_j$ to be justified only by means of epistemically circular arguments; many of them will be justified by other key components. But, nor is every deep disagreement directly deep.
Suppose $K_{i}\ncong K_{j}$ by virtue of some unique key component, say $k\in K_{i}$, such that 
$k$ does not correspond to any element of $K_j$. 
If $k$ has a non-question-begging justification, that justification must also appeal to a key component, that is some other element of $K_{i}$, say $k^{\prime}$. If 
$k^{\prime}$ has a counterpart in 
$K_{j}$, then the disagreement is not deep. But it would now appear that, contrary to our hypothesis, $K_{i}\cong K_{j}$, since the sole obstacle to 
preservation of meaning can be explained in terms common to both theories.
So, if $k$ is sufficient to establish $K_{i}\ncong K_{j}$, then it must lack an independent justification within $T_i$, meeting the Mutual Circularity condition.
\end{enumerate*}

What are some examples of deep disagreements in mathematics?
Silvia Jonas helpfully catalogues ``some of the most important fundamental disagreements in current foundational mathematics'' \cite[304]{Jonas19}. These include the disagreement between classicists and intuitionists as to the appropriate rules for logic (discussed above as an example of a revolution) and disagreements amongst set theorists over the choice of axioms 
and whether axioms supplementary to ZFC are required. Axioms and rules of logic are clearly key components of any theory that depends upon them.
Although Jonas's concern is not with deep disagreement as such, each of these disputes may plausibly be seen as an example thereof. Certainly it is uncontroversial that the competing key components---logical principles or axioms---are epistemic principles. 
However, all of these disagreements occur in the foundations of mathematics, which is remote from the working practice of most mathematicians. Do deep disagreements arise elsewhere in mathematics? To address this question, I will turn to an extended example.


\section{
Inter-Universal Teichm\"uller Theory}

Inter-universal Teichm\"uller theory (IUT or IUTeich, for short) is the work of the 
eminent mathematician Shinichi Mochizuki and has recently been the focus of substantial attention in the mathematical world and beyond. 
The conjecture which provoked 
it, the $abc$ conjecture, is relatively simple to state:
\begin{conjecture}
[
Oesterl\'e--Masser, 1985
]
For every $\varepsilon > 0$, there are only finitely many triples $(a, b, c)$ of co-prime positive integers where $a+b = c$, such that $c > d^{1+\varepsilon}$, where $d$ denotes the product of the distinct prime factors of $abc$.
\end{conjecture}
Suppose you have a set of three natural numbers $\{a,b,c\}$, with no prime factors in common (hence 
they are co-prime), such that the third is the sum of the first two. Generally it is not the case that the third, the sum, will be greater than $d$, 
the product of their distinct prime factors. 
As an arbitrary example, 
take $a$ and $b$ as 15 and 28, and therefore $c$ as 43. They are co-prime, but $d=2\times3\times5\times7\times43=9030$, 
which is, of course, a great deal larger than $43$. 
So those numbers do not work. But 
if you take, say, 1 and 63, their sum is 64 but the product of 
the prime factors of all three numbers is 42, which is less than 64. 
That was a carefully chosen case:
the $abc$ conjecture states that such cases are 
not just scarce, they are finite in number. 
Despite the simplicity of this statement, it resists simple proof.
It is also of wide-ranging and profound significance for very fundamental aspects of mathematics. As the number theorist Michel Waldschmidt puts it, the conjecture
``describes a kind of balance or tension between addition and multiplication, formalizing the observation that when two numbers $a$ and $b$ are divisible by large powers of small primes, $a + b$ tends to be divisible by small powers of large primes'' \cite[quoted in][]{Fesenko16}.
As such, it has many significant consequences---not least Fermat's Last Theorem \cite[1227]{Granville02}.

\subsection{The news from Kyoto}

\begin{table}[tp]
\caption{Timeline of the path to publication of Mochizuki's IUT papers}
\begin{center}
\begin{tabular}{lp{3in}}
\toprule
August 2012&
IUT preprints appear and are submitted to ``a certain mathematics journal'' \cite{Mochizuki20a}.\\
December 2013 \& 2014&Mochizuki posts end of year progress reports on the verification process \cite{Mochizuki13,Mochizuki14}.\\
December 7--11, 2015&
Oxford workshop on IUT. A report from Brian Conrad appears shortly afterwards \cite{Conrad15}.\\

May 2016&
Journal responds: IUT papers accepted with minor corrections; corrected mss returned shortly afterwards \cite{Mochizuki20a}.\\

September 2017&
Journal requests final production corrections \cite{Mochizuki20a}. \\


December 16, 2017&
Japanese press reports that IUT papers are to be accepted at \emph{PRIMS}: ``a decision could be made to publish them \dots\ as early as January'' \cite{Ishikura17}. 
Widespread online discussion ensues, mostly assuming IUT papers have been accepted at \emph{PRIMS} \cite[e.g.][]{Woit17}. 
\emph{PRIMS} issues statement that IUT papers ``have not yet been accepted in a journal''.\\

December 21, 2017&
Peter Scholze first(?) publicly suggests that the proof of Corollary 3.12 contains an inferential gap, in a comment on Frank Calegari's blog \cite{Calegari17}. 
Brian Conrad replies that his December 2015 account of the Oxford IUT meeting provoked ``three independent unsolicited emails'', all indicating a gap at Corollary 3.12.\\

March 15--20, 2018&
Discussion, primarily of Corollary 3.12, between Mochizuki, Yuichiro Hoshi, 
Scholze, and Jakob Stix takes place in Kyoto.\\

May/August 2018&
Scholze and Stix's ``Why $abc$ is still a conjecture" appears \cite{Scholze18}, 
together with responses from Mochizuki \cite{Mochizuki18,Mochizuki18b}. \\

January 5, 2020&
Mochizuki publishes ``whistleblower'' blog post complaining of ``a miserable black hole emerging in the mathematical world" \cite{Mochizuki20a}.\\ 

February 5, 2020&
Acceptance of IUT papers in \emph{PRIMS} \cite[4]{Mochizuki21}.\\

April 3, 2020&
Acceptance 
publicly announced.\\

March 5, 2021&
IUT papers appear in \emph{PRIMS}, with a brief preface.\\
\bottomrule
\end{tabular}
\end{center}
\label{tab:timeline}
\end{table}%

In August 2012, Mochizuki 
uploaded a 500 page sequence of four papers to his website, entitled Inter-Universal Teichm\"uller Theory I through IV. 
These papers contain (what purports to be) a proof of the $abc$ conjecture.
Mochizuki is a highly respected mathematician with a track record of major results, so although he did not 
publicize the papers, they swiftly attracted attention within the mathematical world. 
Further investigation posed a severe challenge, not just because of the papers' length and dependence on an even greater volume of earlier work, but also because of the intrinsic difficulty of the material. 
The initial impression of the mathematical community was bafflement. It took a long time before any sort of consensus began to emerge---if it has. 
Table~\ref{tab:timeline} provides a timeline of some of the main incidents in the IUT papers' path to publication.
By early 2017, 
the process of verification of IUT had been ongoing for about five years. Mochizuki had succeeded in convincing a small circle of other mathematicians, some of them colleagues of his in Kyoto, but also including a number of international visitors. Several conferences had been organized in an attempt to disseminate understanding of IUT more widely, with mixed results.
Although the papers remained unpublished, it was understood that they were going through peer review, a process that can take many years even in much less exceptional cases.
One mathematician I spoke to at this time expressed concern that the results were ``losing freshness'', and more overt scepticism may well have circulated privately, but most public comment was circumspect.
This began to change in December 2017, initially with what seemed like a very favourable development. 

A rumour, 
apparently arising from 
a story in a Japanese newspaper, rapidly spread that Mochizuki's papers had been accepted at \emph{Publications of the Research Institute for Mathematical Sciences} (\emph{PRIMS}). However, this development was seen as controversial for two reasons. Firstly, Mochizuki is the Editor-in-Chief 
of \emph{PRIMS}, which is the in-house journal of the Research Institute for Mathematical Sciences at Kyoto, where he has spent most of his career. 
Nonetheless, \emph{PRIMS} is an eminently respectable journal, it is not unknown for major results to be published in local journals,%
\footnote{For example, the two papers comprising the computer-assisted proof of the Four Colour Theorem were published in the \emph{Illinois Journal of Mathematics}, apparently ``to ensure expert refereeing'' \cite[38]{MacKenzie99}. 
Indeed, Mochizuki has offered a similar rationale, that \emph{PRIMS} is ``by far the most [and indeed perhaps the only truly] \emph{technically qualified}'' journal \cite[4, parenthesis and emphasis in original]{Mochizuki21}.
Of course, when the definition of expert is itself at stake, such arguments 
may appear double-edged: see \S\ref{sec:IUTdeep} below for further discussion.}
and there are procedures by which editors may publish in their own journals while allaying suspicions of impropriety.
\footnote{Procedures which were followed in this case:
``the Editorial Board \dots\ 
form[ed] a special committee to handle it, excluding him [Mochizuki] and with an Editor-in-Chief substituting for him'' \cite[1]{RIMS2021}.}
The second concern was rather more serious: several mathematicians were aware of criticisms that had been communicated to Mochizuki and to which the most recently available revisions of his paper made no apparent reply. Specifically, more than one mathematician independently reported finding a gap in the proof of Corollary 3.12 of the third paper. 
One such person was Peter Scholze, a significantly younger, but also very distinguished mathematician, indeed a 2018 Fields medalist. 
Corollary 3.12 
is several hundred pages in, but 
most of those pages 
are spent introducing a 
complex network of novel 
concepts and setting things up for subsequent theorems which are then stated and given very brief proofs appealing to this vast apparatus of definitions. 
Corollary 3.12 links the results gleaned from this formidable endeavour back to more familiar mathematics. It sets up the fourth paper in which the $abc$ conjecture itself is derived from this result.
Scholze professed that there was an unbridgeable gap at Corollary 3.12 and Scholze's subsequent 
interactions with Mochizuki have not convinced him otherwise. 

It is important to clarify what is meant by the accusation that Mochizuki's proof of Corollary 3.12 contains a gap. Most proofs contain gaps of some sort and most of those gaps are harmless and understood as such by the mathematical community.
Don Fallis draws a very useful distinction between unintentional inferential gaps and intentional enthymematic or untraversed gaps.
Intentional gaps are generally innocuous, although edge cases can 
be controversial.
Enthymematic gaps arise when
the mathematician has worked through the proof in detail, but chosen to omit some steps from the presentation, 
presuming that the missing steps will be unnecessary for the intended audience. 
Such ``skipping'' \cite[259]{Davis72} is a presentational choice. 
Of course, it would be a bad choice if 
the audience find the proof impossible to follow. This sort of misjudgment may be more likely in cutting edge work (such as the 
examples in \S\ref{sec:Comparison} below).
Untraversed gaps occur when a mathematician does not attempt to verify every step of the proof. In some cases, this can be consistent with the acceptance of the proof by the mathematical community. In exceptional cases, even ``universally untraversed gaps'', unverified by \emph{any} mathematician have been admitted in proofs \cite[59]{Fallis03}. 
Subsequent interview-based research with mathematicians bears out Fallis's impression: ``universally untraversed gaps in published proofs are, while not necessarily few, quite innocent'' \cite[246]{Andersen20}.
However, an inferential gap is always a mistake. It occurs when
``the mathematical proposition that the mathematician was trying to prove does not follow by basic mathematical inferences in the manner that the mathematician had in mind'' \cite[51]{Fallis03}.
This is the sort of gap which Scholze diagnosed in the proof of Corollary 3.12.

Accompanied by a colleague, Jakob Stix, Scholze travelled 
to Kyoto in 2018 in an attempt to 
resolve the concerns with Corollary 3.12 in conversation with Mochizuki and his colleagues.
%
%
However, 
Scholze and Stix came away convinced that there was no proof: ``the suggested proof has a problem \dots\ so severe that in our opinion small modifications will not rescue the proof strategy'', as they write in their report of the meeting \cite[1]{Scholze18}. 
This relatively brief manuscript contains a fairly devastating indictment.
Mochizuki's replies 
are much lengthier, but have not convinced either Scholze and Stix or much of the rest of the international mathematical community \cite{Mochizuki18,Mochizuki18b,Mochizuki21}.
Scholze and Stix do stress that they are making ``certain radical simplifications, and it might be argued that such simplifications strip away all the interesting mathematics that forms the core of Mochizuki’s proof'' \cite[4]{Scholze18}.
Their approach is to try and find a way of reconstructing 
Mochizuki's reasoning that does without 
much of his taxonomy: a simplification which 
Scholze and Stix claim (but Mochizuki denies) captures all the key points of Mochizuki's argument.
Ultimately they produce a diagram of relations between various entities and assert that this must be inconsistent, which thereby vitiates Mochizuki's whole project:
``Mochizuki wanted to introduce scalars of $j^2$ somewhere on the left part of this diagram \dots\ 
However, it is clear that {this will result in the whole diagram} having monodromy $j^2$, i.e., {being inconsistent}'' \cite[10]{Scholze18}.
At that point things collapse---at least on Scholze and Stix's version of Mochizuki's argument.

In his reply, Mochizuki strongly rejects that conclusion:
\begin{quote}
In some sense, the \emph{main assertion} of SS underlying this argument in \S2.2 concerning identifications of copies of $\mathbb{R}$ is the following:
\begin{enumerate}
\item[(Lin)] the relationship between any two of these copies of $\mathbb{R}$ is a simple, straightforward \textbf{linear relationship}, given by \emph{multiplication by some scalar}, i.e., \emph{multiplication by some positive real number}.
\end{enumerate}
Here, it should be stated clearly that this assertion (Lin), which underlies the argument of \S2.2, is {\textbf{completely false}} \cite[4; all emphases original in this and subsequent quotations from Mochizuki]{Mochizuki18b}.
\end{quote}
So Mochizuki 
strongly denies 
that Scholze and Stix have successfully linked their simplified version of his argument to his work. As he puts it elsewhere:
\begin{quote}
at numerous points in the March discussions, I was often tempted to issue a response of the following form to various assertions of SS (but typically refrained from doing so!):
\emph{Yes! Yes! Of course, I completely agree that the theory that you are discussing is completely absurd and meaningless, but} \textbf{that} \emph{theory is completely different from IUTch!} \cite[39]{Mochizuki18}.
\end{quote}
Mochizuki later complains that ``adherents of the RCS'', that is the Redundant Copies School, or mathematicians who accept the simplification of Scholze and Stix, are criticizing 
``\emph{logically unrelated fabricated versions} of inter-universal Teichm\"uller theory in which the crucial \textbf{logical AND “$\wedge$”} relation satisfied by the $\mathrm\Theta$-link of inter-universal Teichm\"uller theory is replaced by a \textbf{logical OR “$\vee$”} relation or, alternatively, by a \textbf{logical XOR “$\dot\vee$”} relation'' \cite[23]{Mochizuki21}.

Scholze remained unpersuaded. 
Although he and Stix did not produce a further 
manuscript responding to Mochizuki's reply, he has made more informal public comments expanding on his criticisms, such as this:
\begin{quote}
I may have not expressed this clearly enough in my manuscript with Stix, but {there is just no way that anything like what Mochizuki does can work}. (I would not make this claim as strong as I am making it if I had not discussed this \dots\ 
with Mochizuki in Kyoto for a whole week; the following point is extremely basic, and Mochizuki could not convince me that one dot of it is misguided, during that whole week.) \dots\ 
The reason it cannot work is {a theorem of Mochizuki himself}. This states that a hyperbolic curve $X$ over a $p$-adic field $K$ (maybe with some assumptions, all of which are always satisfied in all cases relevant to IUT) is determined up to isomorphism by its fundamental group $\pi_1(X)$, and in fact automorphisms of $X$ are bijective with outer automorphisms of $\pi_1(X)$. Thus, {the data of $X$ is completely equivalent to the data of $\pi_1(X)$ as a profinite group up to conjugation}. In IUT, Mochizuki always considers the latter type of data, but of course up to equivalence of groupoids this makes no difference 
\cite{Scholze20}.
\end{quote}
Thus we arrive at a point where Scholze accuses Mochizuki of not perceiving the implications of his own results and Mochizuki accuses Scholze of failing to distinguish ``and'' from ``or''.
As Wittgenstein comments of what we now call deep disagreements, 
``Where two principles really do meet which cannot be reconciled with one another, then each man declares the other a fool and heretic'' \cite[\S 611]{Wittgenstein72}.

The false alarm in 2017, when it was rumoured that Mochizuki's papers 
had been accepted, did 
correctly reveal that they were under consideration at \emph{PRIMS}. 
In April 2020, their acceptance at \emph{PRIMS} was duly announced, and the four papers appeared in print the following year.
However, it would seem that Scholze's intervention and his subsequent impasse with Mochizuki may have delayed the process. As the \emph{Asahi Shinbun} reported in its account of the papers' acceptance:
\begin{quote}
In late 2017, it appeared that the articles would be published, but mathematicians in the West pointed out what they considered inappropriate leaps in logic in a core portion of the articles.
That led the journal editorial board to continue with their assessment. A number of other outside experts were consulted and it was only in February that Mochizuki’s proof was considered to no longer have any problems \cite{Ishikura20}.
\end{quote}
Understandably, Mochizuki seems to have found this delay frustrating, complaining of ``a miserable black hole emerging in the mathematical world" \cite{Mochizuki20a}.
Nonetheless, the papers were eventually published 
with, as far as I can determine, 
only incidental corrections from how they had appeared on Mochizuki's website. 
In particular, the proof of Corollary 3.12, which Scholze and Stix claim to contain an unbridgeable gap, is broadly unchanged.

Scholze was invited to review Mochizuki's papers for one of the two major venues for post-publication review in mathematics, \emph{zbMath}.%
\footnote{The erstwhile \emph{Zentralblatt für Mathematik und ihre Grenzgebiete}. Its rival is \emph{Mathematical Reviews}. For further details of this history, see \cite[285 f.]{Barany18}.}
He makes it very clear in his review that his opinion has not changed. Referring to \cite{Scholze18}, he asserts that
``The concerns expressed in this manuscript 
have not been addressed in the published version'' \cite{Scholze21}.
As he summarizes the problem,
``at some point in the proof of Corollary 3.12, things are so obfuscated that it is completely unclear whether some object refers to the $q$-values or the $\mathrm\Theta$-values, as it is somehow claimed to be definitionally equal to both of them, up to some blurring of course, and hence you get the desired result'' \cite{Scholze21}.
Hence
``the argument given for Corollary 3.12 is not a proof, and the theory built in these papers is clearly insufficient to prove the ABC conjecture'' \cite{Scholze21}.
The other major review journal, \emph{Mathematical Reviews}, commissioned their review of the IUT papers from Mohamed Sa\"idi, an ally of Mochizuki's and frequent visitor to Kyoto \cite{Saidi22}. He makes no mention of Scholze's concerns.

So the mathematical world seems to be in what the American number theorist Frank Calegari anticipated as
``the ridiculous situation where ABC is a theorem in Kyoto but a conjecture everywhere else'' \cite{Calegari17}.
Certainly, the news from Kyoto is that this is now unequivocally a theorem. 
Popular treatments of Mochizuki's work published in Japanese 
bear out this interpretation. 
For example, even before the official acceptance of Mochizuki's papers, Fumiharu Kato, 
a mathematics professor 
and popular writer on mathematics, 
published a book in Japanese with the English title: \emph{Mathematics that Bridges Universes: The Shock of IUT Theory}. This promises, optimistically one might think,
``that {IUT theory is rooted in a natural way of thinking} that can be understood by ordinary people who are not mathematicians'' \cite{Kato19}.
Another popular treatment, this time of just the $abc$ conjecture rather than Mochizuki's results, 
although clearly inspired by that work, refers to it as 
``{a super-difficult problem in mathematics whose proof has been confirmed after the 2020 peer review}
'' \cite{Koyama21}.
So, as far as the Japanese public is concerned, it would seem that the $abc$ conjecture has been proved by their man. 
Indeed, Scholze's comments 
bemoan that his disagreement with Mochizuki seems to have become a matter of national chauvinism: 
``I’m really frustrated with the current situation \dots\ 
effectively arguing along national lines; again, this strikes deep into my heart. I’m really quite surprised by the strong backing that Mochizuki gets from the many eminent people (who I highly respect) at RIMS'' \cite{Scholze20}.

Mochizuki does have some defenders outside Japan. One of the most outspoken is 
Ivan Fesenko, a mathematics professor 
at the University of Nottingham in the UK, who 
has responded to the criticisms of Scholze and others with some asperity:
\begin{quote}
Some, affected by negative emotions, broke professional rules of conduct and made public their ignorant and sometimes intolerant opinions. Tellingly, the only questions produced were shallow and misplaced \dots\ 
It should be further analysed why has it become acceptable for some mathematicians to produce public negative statements about the work of other mathematicians without any concrete mathematical justification and what fundamental harm to pioneering research this unethical behaviour makes \cite[5 f.]{Fesenko18}. 
\end{quote}
Fesenko is also one of the very small number of mathematicians to have published a survey article that attempts to clarify 
Mochizuki's results \cite{Fesenko15}.
However, this paper does not claim to be a full recapitulation of the IUT papers, and the consensus of the mathematical community is that much more work in this vein is required.
For example, Brian Conrad 
contrasts Mochizuki's work with that of Scholze on perfectoid spaces:
``the production of efficient survey articles and lectures getting right to the point in a moderate amount of space and time occurred very soon after that work was announced. 
Everything I understood during the week in Oxford supports the widespread belief that there is no reason the same cannot be done for IUT, exactly as for other prior great breakthroughs in mathematics'' \cite{Conrad15}.
While Conrad's optimism now seems premature, 
there is ongoing work of the kind he suggests: perhaps the most promising is the sequence of papers by Taylor Dupuy and Anton Hilado \cite[beginning with][]{Dupuy20a}.

\subsection{Comparison with other results}\label{sec:Comparison}

Several aspects of the reception of Mochizuki's IUT papers as outlined above may initially arouse suspicion. 
In particular, that the reasoning is found hard to follow;
that the proof is believed by at least some of its critics to contain a gap;
and that widespread understanding of the result requires a comprehensive rewrite by later, independent mathematicians.
However, while these are perhaps uncommon features amongst mathematical proofs in general, they
echo the vicissitudes of many other major mathematical results. (See Table~\ref{tab:comparison}.)
For example, Kurt Heegner's (1952) proof that there are exactly nine complex quadratic fields of class-number one has been described as ``written in the most horrible style that you can think of'' \cite[40]{Kolata83};
was presumed to inherit an inferential gap from 
earlier work by Heinrich Weber known to contain errors;
and the result was therefore reproved independently by other mathematicians.
But one of those mathematicians then reexamined Heegner's proof and concluded that the dependence on Weber could be easily avoided and  
``that there is in fact only a very minor gap in Heegner’s proof'' \cite[16]{Stark69}. 
A much more prominent result with similar parallels is Michael Freedman’s \citeyearpar{Freedman82} proof of the four-dimensional Poincaré conjecture. 
In stark contrast with the protracted reviewing process of the IUT papers, this proof was published with only minimal peer review.
After Freedman presented his proof at a conference, the paper was solicited by a journal which sent it to a close associate of Freedman’s; when the reviewer vouched for the proof's correctness but complained that a full report would take many weeks of work (``each page in the paper 
\dots\ generated at least a page of notes, questions, corrections, and typos''), the journal accepted the paper without corrections \cite{Kirby20}.
Perhaps as a result of this unorthodox path to publication, the paper has long been seen as particularly difficult to follow, although Freedman’s personal explanations are apparently much clearer (which may be why his proof has been generally accepted as correct) \cite{Hartnett21}. 
This asymmetry eventually provoked the publication of a lengthy, multiply authored book recasting Freedman’s arguments in accessible terms \cite{Behrens21}.
Louis de Branges’s (1984) proof of the Bieberbach conjecture exhibits a similar pattern: 
\begin{quote}
[De Branges] completed a manuscript of 385 pages \dots\ 
As he tells it, he was disappointed that the U.S. mathematicians to whom he had sent his manuscript had not yet been able to verify his long proof. In Leningrad, de Branges presented his work to the members of the seminars in functional analysis and geometric function theory. In a large number of sessions, the proof was verified and some inessential errors corrected. Finally, through hard work under de Branges' direction, a relatively short proof of the Lebedev-Milin conjecture [the crucial step required to prove the Bieberbach conjecture] was distilled from the original manuscript \cite[513]{Korevaar86}.
\end{quote}
Several subsequent authors reworked de Branges’s proof into even shorter and simpler forms, culminating in ``a wallet-sized, high-school-level proof'' \cite[165]{Krantz11}.
Most conspicuous of all (but also perhaps most idiosyncratic) is
Grigori Perelman’s (2002) proof of the three-dimensional Poincaré conjecture.
Perelman never submitted his preprints for publication, but they were independently vetted by two teams of mathematicians. Nonetheless, the Chinese mathematician Shing-Tung Yau controversially claimed that the proof was only completed in an article reconstructing Perelman’s work (and written by two of Yau's prot\'eg\'es) 
\cite{Nasar06}.
A subsequent correction to this article withdrew the claim \cite[244]{Szpiro07}.
As we have seen, Mochizuki's IUT papers have caused widespread bafflement and have been explicitly accused of harbouring an inferential gap. Several mathematicians have expressed the hope that the core ideas will eventually be rewritten by others in a more accessible form, a project which is still in its very early stages. 
The point of this 
series of comparisons is to stress that none of these factors should necessarily count against Mochizuki's work. Of course, neither are they evidence in its support.

\begin{table}[tp]
\caption{Comparison of the reception of several major mathematical results}
\begin{center}
\begin{tabularx}{.95\textwidth}{l*{5}{Y}}
\toprule
&Resists&Suspected&Rewritten&Deep\\
&understanding&of gap(s)&by others&disagreement&Revolution\\
\midrule
Heegner (1952)&\checkmark&\checkmark&\checkmark&$\times$&$\times$\\
Freedman (1982)&
\checkmark&$\times$&\checkmark&$\times$&{\checkmark}\\
De Branges (1984)&
\checkmark&?&\checkmark&$\times$&?\\
Perelman (2002)&\checkmark&\checkmark&\checkmark&$\times$&{\checkmark}\\
Mochizuki (2012)&\checkmark&\checkmark&
ongoing&\checkmark&\checkmark\\
\bottomrule
\end{tabularx}
\end{center}
\label{tab:comparison}
\end{table}%

Which, if any, of these earlier results are 
revolutionary or give rise to deep disagreements?
Freedman and Perelman's work has often been characterized as revolutionary and is marked by substantial conceptual innovations which ensure that it meets the narrower definition set out in \S\ref{sec:Revolutions} above.
Heegner, on the other hand, was working within an established tradition. The same can probably be said of de Branges's proof, despite its higher profile.
By contrast, few if any of these cases could be said to give rise to a deep disagreement.
Heegner died in obscurity before his work was rediscovered, so there was little opportunity for disagreement. Moreover, the 
misperception of his proof as flawed turned on a relatively superficial confusion.
Freedman's proof, although notoriously hard to follow, has generally been accepted as correct since its inception.
De Branges's original manuscript seems to have been ignored rather than critiqued; if its recipients suspected it of containing gaps, they may have done so 
without reading it, making only for a very superficial disagreement.
Yau certainly accused Perelman's preprints of containing substantial gaps, which entails some sort of disagreement, even if the consensus of the mathematical community has been against him. 
However, Yau accepted Perelman's proof strategy and was not claiming that the alleged gaps were irreparable---on the contrary, he asserted that they had been subsequently filled. Essentially this  disagreement concerned the level of credit due the mathematicians who filled these gaps: that is a priority dispute, not a deep disagreement.

\subsection{IUT as revolutionary}\label{sec:IUTrev}
We are now in a position to assess IUT's claims to revolutionary status.
Some of Mochizuki's defenders, such as
Fesenko, 
have overtly (and repeatedly) ascribed a revolutionary nature to 
his work:
``This theory is so radically different from anything that came before it that it is natural to ask whether it will induce a paradigm shift'' \cite[436]{Fesenko15};
``We had mathematics before Mochizuki’s work---and now we have mathematics after Mochizuki’s work'' \cite[quoted in][181]{Castelvecchi15};
``IUT is different in its philosophy and main ideas from anything we have known in conventional number theory. It is already changing mathematics, and as more people learn and develop IUT, this will continue'' \cite{Fesenko16}. 
But most importantly, Fesenko has buttressed these claims by appeal 
to IUT's taxonomic incommensurability with other mathematical theories: it is 
``a highly novel exotic theory with a two-digit number of new concepts'' \cite{Fesenko18}.
Other mathematicians have also ascribed ``revolutionary new ideas'' to IUT \cite[Jeffrey Lagarias, cited in][14]{Castelvecchi16}.
If this characterization is correct, IUT would meet the definition of (paraglorious) revolution proposed in 
\S\ref{sec:Revolutions} above.

Fesenko has also stated 
that Mochizuki compares himself to Alexander Grothendieck, a major twentieth-century mathematician whose work is widely seen as revolutionary
\cite[181]{Castelvecchi15}. 
This may be an immodest claim, but it is an understandable one: Grothendieck's work is an obvious and important influence on Mochizuki. His first major result proved an open conjecture made by Grothendieck, 
his research is largely within fields which Grothendieck worked in (or indeed invented),
and Mochizuki's expository style, at least within the IUT papers, clearly echoes that of Grothendieck.
The hundreds of pages of definitions and theorems with brief proofs leading up to the fateful Corollary 3.12 at least stylistically resemble Pierre Deligne's description of 
``a typical Grothendieck proof as a long series of trivial steps where `nothing seems to happen, and yet at the end a highly non-trivial theorem is there'\,'' \cite[quoted in][302]{McLarty03}.
On the other hand, for others,
\begin{quote}
the comparison between Mochizuki and Grothendieck is a poor one. Yes, the Grothendieck revolution upended mathematics during the 1960's `from the ground up'. But \dots\ 
the success of the Grothendieck school is not measured in the theorems coming out of IHES in the '60s but in how the ideas completely changed how everyone in the subject (and surrounding subjects) thought about algebraic geometry \cite{Calegari17}.
\end{quote}
Mochizuki's ideas have not yet had anything like this sort of influence. Even if we take the peer-reviewed publication of his papers 
as establishing the cogency of his arguments (an issue I return to below), the true revolutionary potential of his project would turn on the wider application of his techniques. While his supporters insist that this will follow, they have few such examples so far, even after nearly a decade.

Both of these arguments for the revolutionary character of IUT presume its success.
If it does succeed, Mochizuki's approach 
exhibits taxonomic incommensurability with 
prior work in algebraic geometry, and is thereby paragloriously revolutionary.
Even if IUT 
is eventually established beyond dispute to be fundamentally inconsistent, 
I have argued elsewhere that 
the inter-universal Teichm\"uller theory paradigm would, like other failed architectural conjectures, undergo a form of inglorious revolution as it collapses \cite[147]{Aberdein17}.
So there is some sort of revolution involved one way or the other: 
a paraglorious revolution if it succeeds; or an inglorious revolution if it fails.
Which of these situations obtains remains, of course, the subject of a substantial 
disagreement, which I address 
next.

\subsection{IUT and deep disagreement}\label{sec:IUTdeep}
What are Mochizuki and Scholze disagreeing about and is their disagreement deep? Most prominently, but also superficially, they disagree as to whether the IUT papers comprise a proof of the $abc$ conjecture.
But, as we have seen, the disagreement is more focussed: does Corollary 3.12 contain an inferential gap?
That disagreement turns on the legitimacy of Scholze's so-called ``redundant copies'' reconstruction of Mochizuki's work (RCS-IUT, as the latter calls it). Mochizuki makes clear that he agrees with Scholze that 
``RCS-IUT is indeed a meaningless and absurd theory that leads immediately to a contradiction'' \cite[5]{Mochizuki21}.
For Mochizuki, the redundant copies are not redundant: rather it is their subtle interplay on which his whole proof relies.
As the mathematician David Roberts 
cautiously observes,
\begin{quote}
It remains entirely possible that those radical simplifications engineered by Scholze and Stix identified objects that are isomorphic only after \emph{some} stage of a tower of forgetful functors,%
\footnote{In category theory, a forgetful functor maps a more complex 
to a less complex algebraic structure, by ``forgetting'' some of that complexity.}
but not at the earlier stage at which they were meant to be considered. A system of objects may have been identified with a different, simpler system of objects \emph{unnaturally}, various necessary compatibility conditions violated \cite{Roberts19}.
\end{quote}
That is what Mochizuki insists upon and what Scholze denies. Hence the substance of their dispute concerns the exact nature of the relation between IUT and RCS-IUT: are they interchangeable, as Scholze asserts, or fundamentally distinct, as Mochizuki would have it?
Of course, even if this dispute were settled in Mochizuki's favour, his proof of the $abc$ conjecture might still fail for other reasons. As Roberts continues,
\begin{quote}
Even if Scholze and Stix’s analysis is flawed, and Mochizuki’s categorical foibles are harmless, his papers may still have a gap, some innocuous assumption unchecked, some existence statement unjustified---an $abc$-sized gap deep in the proof of Theorem 3.11 or Corollary 3.12. This will only be found by careful study and ideally a rewriting of Mochizuki’s papers into more standard language 
\cite{Roberts19}.
\end{quote}
The discovery of such a gap would render the Mochizuki--Scholze disagreement moot.
But, so far, neither party has produced anything likely to be seen as a watertight proof of either side of the dilemma, hence the disagreement endures.

The relationship between IUT and RCS-IUT is ultimately a matter of fact, rather than epistemic principle.
However, as Michael Lynch observes,
``debates over the basic facts can themselves turn into disagreements over whose sources and standards for facts are trustworthy'' \cite[149]{Lynch20}.
Just such a shift has occurred in this disagreement, at least on Mochizuki's side.
He attributes RCS-IUT (and what he perceives as its misleading conflation with IUT) to
``a growing collection of mathematicians who have a somewhat \emph{inaccurate} and \emph{incomplete}---and indeed often quite \emph{superficial---understanding} of certain aspects of the theory
'' \cite[14]{Mochizuki21}. Mochizuki asserts that
``the number of professional mathematicians who have achieved a sufficiently detailed understanding of inter-universal Teichm\"uller theory to make independent, well informed, definitive statements concerning the theory \dots\ 
is roughly on the order of 10'' \cite[11]{Mochizuki21}.
This in turn justifies his choice of a local journal as ``by far the most [and indeed perhaps the only truly] \emph{technically qualified}'' venue to referee the IUT papers \cite[4]{Mochizuki21}.
Scholze and other sceptics of IUT 
dispute the implication that anyone outside Mochizuki's team of ten lacks sufficient expertise to judge the papers.

As we saw in \S\ref{sec:IUTrev}, IUT is widely seen as taxonomically incommensurable with mainstream mathematics, which would ensure not only that it is revolutionary, but also that a disagreement over its success meets the definition of deep disagreement proposed in \S\ref{sec:DD}.
We can also see that the Mochizuki--Scholze dispute meets Lynch's criteria for deep disagreement.
The two parties have a common epistemic goal, to determine whether the IUT papers constitute a proof of the $abc$ conjecture, thereby satisfying Commonality.
The two parties affirm distinct epistemic principles, at least in the sense that they disagree over who is an expert, which satisfies Competition.
Non-arbitration requires further discussion: there are two candidate arbiters, but neither of them are helpful here. Firstly, 
proof is often taken to be reducible to formal derivation, which is susceptible to purely mechanical (or even automated) checking. This would seem to render any apparent disagreements over what counts as a proof easily resolvable. 
However, this account is misleading: in practice few mathematical proofs have ever 
been successfully represented as formal derivations. 
Doing so requires an exceptionally fine-grained grasp of every detail of the informal proof.
Resolving the IUT/RCS-IUT dilemma would thus be an essential step in any such formalization of Mochizuki's work. 
Formalization cannot be expected 
to arbitrate a dispute that would have to be settled before the formalization could begin.
A second, more
conventional, yardstick for settling disagreements about the status of proofs is peer review. The IUT papers have now been published in a peer-reviewed journal. However, this has not settled the disagreement either.
As Calegari observed before the papers were published,
``whether the papers are accepted or not in a journal is pretty much irrelevant here; it’s not good enough for people to attest that they have read the argument and it is fine, someone has to be able to explain it'' \cite{Calegari17}.
Calegari implicitly endorses the useful maxim that
``Proof $=$ Guarantee $+$ Explanation'' \cite{Robinson00,Brown17a}. Even if we were to grant that the peer review process has provided the first component (which in turn would require endorsing something like Mochizuki's 
view on the scarcity of expertise in IUT),
we would still be lacking the second component:
\begin{quote}
[The referee process] forces the author(s) to bring the clarity of the writing up to a reasonable standard for professionals to read it (so they don’t need to take the same time duration that was required for the referees, amongst other things). This latter aspect has been a {complete failure}, calling into question both the quality of the referee work that was done and the judgement of the editorial board at PRIMS to permit papers in such {an unacceptable and widely recognized state of opaqueness} to be published 
\cite{Calegari17}.
\end{quote}
As Calegari suggests, if the referees have not achieved one of the components, that in turn undermines our trust that they have achieved either component.
Thus peer review has not resolved the dispute either. In the absence of any other candidate arbiters, Lynch's Non-arbitration criterion has been met.
Lastly, Mutual Circularity also appears to obtain. 
Mochizuki explicitly appeals to his own expertise and that of his closest colleagues 
in support of his proof; much of the rest of the mathematical world appeals to Scholze's acknowledged expertise in the wider context of number theory in support of his claim to have found a gap. 
If the dispute concerns who is a suitably credentialed expert to judge the soundness of the IUT papers, such appeals to expertise must appear circular.



\section{A Virtuous 
Approach to Vicious Disagreements}
I will conclude by turning to the intellectual virtues. 
I wish to show they can help to provide a 
framework for understanding deep disagreements such as that between Mochizuki and Scholze. 
In this task I am building on two prior bodies of work.
Firstly, there has recently been a growing interest in the application of virtue epistemology to questions in the philosophy of mathematics \cite[as surveyed in][]{Aberdein21b}.
Secondly, in other work I have 
defended a virtue approach to argument and, more narrowly, to deep disagreement.
In particular, I have focussed on the vice of 
arrogance 
and the virtue of 
courage 
\cite{Aberdein20a,Aberdein19b}. 
Although other virtues and vices could profitably be discussed in this context, this choice was not 
arbitrary. 
Arrogance and courage may be contrasted as the vice and virtue characteristic of the ``steadfast'' view in the epistemology of disagreement; conversely, its polar opposite, the ``equal weight'' view, may be thought to exemplify either humility or cowardice, depending on one's perspective \cite[\S5.2]{Frances18}.
If, when I learn of a disagreement with my view, I respond by sticking to my guns, I exhibit courage but risk arrogance. Conversely, if I respond to the disagreement by 
weakening my support for my own view, I demonstrate humility but may be succumbing to cowardice.

As I have argued elsewhere, it is helpful in thinking about the role of courage in disagreement to
distinguish 
two species of courage \cite[1210]{Aberdein19b}.
For Per Bauhn ``the courage of creativity (determination in the face of adversities, resourcefulness under pressure) supports the courage of conviction (a sense of responsibility based on moral beliefs), while, simultaneously, creativity feeds on conviction'' \cite[136]{Bauhn03}.
Both involve overcoming fear, whether the fear of failure 
 in the former case or the fear of harm, including social harm, in the latter.
 This distinction is often framed as one between intellectual and moral courage.
Virtue epistemologists often define intellectual courage in terms of intellectual perseverance \cite[689]{Battaly17b}.
For example, for Jason Baehr, intellectual courage is 
``{a disposition to persist} in or with a state or course of action aimed at an epistemically good end despite the fact that doing so involves an apparent {threat to one's own well-being}'' \cite[177]{Baehr11a}.
Moral courage, on the other hand, involves overcoming some sort of danger, specifically where the danger takes the form of 
%
``{a threat to one's social standing}, financial prospects, {relations with one's colleagues}, approval of one's constituents, and so forth'' \cite[107]{Walton86}.
Saliently, Douglas Walton, whose characterization this is, also notes that
``The person of moral courage is open to persuasion and reasonable discussion, but will not give in to pressures until convinced the path is right. {Compromise} is therefore not intrinsically a sign of weakness or cowardice---it could in some cases actually be {a mark of courage}'' \cite[128]{Walton86}.
This quality can be especially valuable in an arguer.
%
Daniel Cohen describes as an ``argument provocateur'' the sort of arguer who is willing to engage in argument to a fault,
starting arguments even with those with whom they agree.
%
However irritating such behaviour can be, it does exemplify 
moral courage in the face of the dangers of argument:
``In some circumstances, arguing is bad form. If we are too sensitive to that, we can become 
(to resort to the notorious, but occasionally apt, war metaphor for arguments) 
\emph{gun shy} about arguing
'' \cite[64]{Cohen05}.
%

Reflection on the Mochizuki--Scholze dispute should convince us that both mathematicians have exhibited courage, although of different kinds.
It should be uncontroversial that Mochizuki has at least shown intellectual 
courage in persevering with 
IUT through its lengthy development and many subsequent tribulations.
It should be uncontroversial that Scholze has at least shown moral courage in acting on his conviction that IUT is vitiated by a fundamental error that the mathematical community risked ignoring. 
While we might hope for consensus in these two judgments, any further assessment will turn on one's perspective on the dispute.
From Mochizuki's perspective, Scholze overestimates his own expertise with respect to IUT, in other words he has not persevered sufficiently with its study, but given up too early and nonetheless concluded that he is sufficiently knowledgable to diagnose irreparable faults. This vice of deficiency with respect to intellectual courage has in turn led him into a vice of excess with respect to moral courage, a morally reckless attack on IUT.
But from Scholze's perspective, Mochizuki's perseverance with IUT is itself a vice of excess with respect to intellectual courage 
and has in turn led him into a vice of deficiency with respect to moral courage, a refusal to admit his mistakes, despite (what Scholze takes to be) overwhelming evidence that IUT does not achieve what Mochizuki had hoped it would.
Thus we arrive at a symmetrical impasse: each disputant sees the other as manifesting 
vices that mirror those the other party is perceiving in him.

What of arrogance and humility?
Helpfully, Colin Rittberg has explored the role of humility in understanding the IUT saga. 
He has discovered yet another symmetrical impasse: 
``Mochizuki’s critics demand that he manifest more intellectual humility by going further in owning the limitations of his communicative efforts to foster understanding of his $abc$-proof. Mochizuki counters by demanding that his critics manifest more intellectual humility in their engagement with his proof'' \cite[5591]{Rittberg21}.
Lynch's work suggests that focussing on arrogance, the corresponding vice, may be more insightful. He notes that
``deep epistemic disagreements may cause the participants to perceive each other more negatively from the epistemic view, including possibly seeing each other as arrogant'' \cite[152]{Lynch20}.
This follows from a specific feature of his account of deep disagreement. 
As we have seen, Mutual Circularity is one of Lynch's hallmarks of deep disagreement: parties to such disagreements end up justifying their epistemic principles through circular reasoning, for want of anything better. 
Such circular arguments may appear innocuous from the perspective of someone who implicitly accepts the epistemic principles in question, and have been defended as such by some epistemologists.
But, from an outside perspective, that of someone who rejects the relevant epistemic principles,  appealing to precisely those principles in their own defence may well appear not only unsound, but arrogantly so.
Something much like the situation Lynch describes seems to have arisen with respect to IUT. As noted above, the dispute between Mochizuki and Scholze exhibits Mutual Circularity, since each party seeks to buttress their position by appeal to expertise, although that expertise is at the root of the deep disagreement. 
Each party thereby argues 
in a fashion that seems quite reasonable from their own perspective, but appears arrogantly circular to the other camp.

One of the limitations of vice epistemology, of which these reflections are a part, is that accusations of epistemic vice can themselves be counterproductive and inflammatory.
On Ian James Kidd's account, 
``the efficacy of a vice-charge is contingent on consensus between critic and target. There must be consensus, first, on the {definition} of the vice being invoked \dots\ 
and, second, on whether the target does in fact {exemplify} that vice'' \cite[192]{Kidd16a}.
When critic and target are in deep disagreement it would appear vanishingly unlikely that they should arrive at such a consensus.
Moreover, vice charging can easily make matters worse, since a failed vice charge ``reflects badly on the critic, especially if they are claiming the moral or epistemic high ground'' \cite[185]{Kidd16a}.
These considerations suggest that vice epistemology will be of limited use to the parties to a deep disagreement, even if it is helpful to philosophers attempting to analyse the dispute.

I have argued elsewhere that humility and courage can help to resolve deep disagreements, or at least, that the corresponding vices of arrogance and cowardice can make such resolutions harder to achieve \cite{Aberdein20a,Aberdein19b}. In part, this is because there are persuasive strategies that can help to resolve deadlocked arguments, but which are open to abuse. Only when such strategies are pursued virtuously can a satisfactory outcome ensue.
However, such an approach may be of less use in the IUT dispute, and in similar mathematical deep disagreements. That is because this case and others like it lack the action-forcing constraints present in many deep disagreements: there is no deadline by which some decision must be taken, some policy adopted, or so forth.
It may be an embarrassment to mathematics that so prominent a peer-reviewed proof is so widely disputed, but it does not become any more of an embarrassment 
by enduring. As such, there is little incentive for mathematicians to seek an immediate resolution.
More technically, the types of dialogue in which the IUT papers and their path to acceptance
belong are inquiry and persuasion dialogues, not negotiation or deliberation dialogues (the distinction is due to \citealp{Walton95}; see \citealp[160 ff., for application to mathematics]{Aberdein21a}).
The latter sort of dialogues seek practical, short-term accommodations, not the sort of stable resolutions to which mathematics aspires.
Hence the proper attitude is one of patience: 
perhaps the ongoing projects to clarify Mochizuki's work or extend it to new results will meet with general approval; perhaps Scholze's critique will be reframed in even more unignorable terms.

\section{Conclusions}
I have proposed criteria for recognising both {revolutions} and {deep disagreements} in mathematics.
We have seen that inter-universal Teichm\"uller theory, Mochizuki's programme, meets the criteria for a mathematical revolution, since it is not taxonomically commensurable with the programme to which it is a successor. The exact nature of that revolution turns on its eventual fate.
It is paragloriously revolutionary if it succeeds 
and ingloriously revolutionary if it fails. 
In either case, the dispute over whether IUT does succeed is a 
deep disagreement. 
It is thereby an example of a deep disagreement arising in mathematical practice rather than in the foundations of mathematics.
I also touched on the relevance of some intellectual virtues to this dispute. 
Intellectual courage
is essential to the defence of unpopular views. 
Defending a position in the face of technical difficulties and adverse attention takes a significant amount of courage. But moral courage 
is also essential for acknowledgement of defeat. 
I have argued elsewhere that disagreement amelioration strategies require close attention to virtues of argument. 
However, such strategies find readiest application when there is a pressing need to find even a temporary solution. Pure mathematics has the luxury of longer horizons and seeks 
permanent outcomes. Although it often maintains a remarkable level of consensus, these priorities do mean that when that consensus breaks down, as it must in a revolution, it may be very slow to heal.

\section*{Acknowledgements}
I am grateful to Colin 
Rittberg and Deniz Sarikaya for the invitation to contribute to this special issue and its associated workshop.
I am further indebted to two reviewers for this journal whose comments 
resulted in some significant improvements.
Fenner Tanswell also provided helpful feedback on an earlier draft.
Much of the research for this article was conducted during sabbatical leave from Florida Institute of Technology. I am thankful 
for this support. 

\end{document}